\documentclass{amsproc}
\usepackage{amssymb,amsthm}
\usepackage{amsmath}
\usepackage{amscd}
\usepackage{subfigure}
\usepackage{epsf,graphicx}
\usepackage[all]{xy}

\newtheorem{thm}{Theorem}
\newtheorem{cor}[thm]{Corollary}
\newtheorem{prop}[thm]{Proposition} 

\newtheorem{defi}[thm]{Definition}

\DeclareFontFamily{OMS}{rsfs}{\skewchar\font'60}
\DeclareFontShape{OMS}{rsfs}{m}{n}{<-5>rsfs5 <5-7>rsfs7 <7->rsfs10
}{} \DeclareSymbolFont{rsfs}{OMS}{rsfs}{m}{n}
\DeclareSymbolFontAlphabet{\scr}{rsfs}

\newcommand{\des}{\displaystyle}

\newcommand{\im}{{\rm{im}}}

\newcommand{\en}{{{\rm End}}}
\newcommand{\ma}{{\rm Map}}
\newcommand{\Co}{{\rm Cob}}
\newcommand{\de}{{\rm dg\mbox{-}vect}}

\newcommand{\Mfu}{{\rm MFunc}}

\newcommand{\Ho}{{\rm H}}
\newcommand{\hg}{{\rm HG}}

\newcommand{\Si}{{\rm S}}
\newcommand{\Ca}{{\rm C}}
\newcommand{\sC}{\scr{C}}

\newcommand{\ob}{{\rm Ob}}

\newcommand{\al}{{\rm A}}

\newcommand{\OB}{{\rm ob}}
\newcommand{\Sc}{{\rm Schur}}
\newcommand{\ve}{{\rm vect}}
\newcommand{\HL}{{\rm HLQFT}}

\newcommand{\I}{{\rm I}}
\newcommand{\id}{{\rm id}}

\renewcommand{\des}{\displaystyle}
\def\N{{\Bbb N}}

\def\R{{\Bbb R}}

\def\R{{\mathbb R}}

\def\C{{\mathcal{C}}}

\def\d{{\mathcal{D}}}
\def\f{{\mathcal{F}}}

\numberwithin{equation}{section}

\begin{document}

\title{Homological matrices}

\author{Edmundo Castillo}
\address{Departamento de Matem\'aticas, Universidad Central de Venezuela (UCV).}
\curraddr{Universidad Central de Venezuela,
Av. Los Ilustres, Los Chaguaramos,
Caracas-Venezuela
}
\email{ecastill@euler.ciens.ucv.ve}

\author{Rafael D\'\i az}
\address{Departamento de Matem\'aticas, Universidad Central de Venezuela (UCV)} \email{rdiaz@euler.ciens.ucv.ve}

\curraddr{Universidad Central de Venezuela, Av. Los Ilustres, Los
Chaguaramos, Caracas-Venezuela }
\date{September 2005}


\begin{abstract}
We define homological matrices, construct examples of
one-dimen\-sion restricted  homological quantum field theories,
and show a relationship between the two theories.
\end{abstract}

\maketitle
\section{Introduction}

The purpose of this paper is to introduce examples of
one-dimensional homological quantum field theories $\HL,$ and to
define a higher dimensional homological analogue of matrices.
First we review the notion of $\HL$ introduced in \cite{Cas} which
is based on the concept of cobordisms dotted with homology classes
of maps into a fixed compact oriented smooth manifold. We demand
that the maps from cobordisms to the fixed manifold be constant on
each boundary component. This restriction is necessary in order to
define composition of morphisms using transversal intersection on
finite dimensional manifolds. We define a HLQFT as a monoidal
functor from that category of extended cobordisms into the
category of vector spaces.

{\noindent}Throughout this paper we work at the chain level. To do
so we introduced in Section \ref{sec1}  transversal $1$-categories
following ideas by     May and Kriz \cite{KM}, and homology with
corners for oriented smooth manifolds.  In Section \ref{sec2} we
defined homological quantum field theory in arbitrary dimensions.
 Section \ref{sec3} contains our main example of one-dimensional
$\HL$ which is based on the notion of parallel transportation
along chains of paths.

{\noindent}In Section \ref{sec4} we introduce a higher dimensional
homological analogue of matrices using the graphical
representation of matrices introduced in \cite{RE}. We describe
the relation between  $1$-dimensional restricted $\HL$ and higher
dimensional matrices.

\section{Transversal $1$-categories}\label{sec1}

Fix a field $k$ of characteristic zero. We denote by $\de$ the
category of differential $\mathbb{Z}$-graded $k$-vector spaces.
Objects in dg-vect are pairs $(V,d),$ where $V=\bigoplus_{i\in
\mathbb{Z}}V_{i}$ is $\mathbb{Z}$-graded $k$-vector space and
$d=\bigoplus _{i\in \mathbb{Z}}d_{i}\colon V \to V$ is such that
$d_{i}\colon V_{i}\to V_{i-1}$ and $d^{2}=0.$ If $V$ is $\de$ then
$V[n]$ is the $\de$ such that $V[n]^{i}=V^{i+n}.$ We say that $V[n]$
is equal to $V$ with degrees shifted down by $n.$ If $v\in V^{i}$
the are write $|v|=i.$ If $V$ is a dg-vect $\Ho(V)$ denotes the
homology of $V.$ A {\it dg-precategory} $\C$ consist of a collection
of objects $\ob(\C)$ and for any $x,y \in \ob(\C),$ a dg-vect
$\C(x,y)$ called the space of morphisms  from $x$ to $y.$ A
dg-category  $\C$ is a dg-precategory together with complex
morphisms $\circ \colon \C(x,y) \bigotimes \C(y,z) \longrightarrow
\C(x,z)$ and $\id_x \colon k \rightarrow  \C(x,x)$ turning $\C$ into
a category. We think of $k$ as the complex with zero differential
concentrated in degree zero. For a dg-precategory $\C,$ we define
the precategory $\Ho(\C)$ as follows: $\ob(\Ho(\C))=\ob(\C),$ and
for any  $x,y \in \ob(\Ho(\C)), \ \Ho(\C)(x,y)=\Ho(\C(x,y)).$ One
checks easily that if $\C$ is a dg-category, then $\Ho(\C)$ is also
a category.

\vspace{0.3cm} {\noindent} A {\it transversality structure}
$\C^{\star}$ in a dg-precategory $\C$ consist of the following data:
\begin{itemize}
\item{For $x_{0},\cdots ,x_{n} \in \ob(\C)$ a $\de$ $\C(x_{0},\cdots
 ,x_{n}).$ }
\item{ inclusions $i_{n}\colon \C(x_{0},\cdots ,x_{n}) \hookrightarrow
\displaystyle\bigotimes_{i=1}^{n}\C(x_{i-1},x_{i}).$ }
\end{itemize}

This data should satisfy the following properties:

\begin{itemize}
\item{$i_{n}\colon \C(x_{0},\cdots ,x_{n}) \to\des
\bigotimes_{i=1}^{n}\C(x_{i-1},x_{i})$  is a quasi-isomorphism,
i.e., the induced map \\  $\Ho({i}_{n})\colon\Ho(\C( x_{0},\cdots
,x_{n})) \to \displaystyle\bigotimes_{i=1}^{n} \Ho(\C(x_{i-1},
x_{i}))$ is an isomorphism.}
\item{ For each partition $n=n_{1}+ \cdots + n_{k},$ the map $i_{n}$
factors through \\
$\displaystyle\bigotimes_{i=1}^{k}\C(x_{n_{i-1}},\cdots
 ,x_{n_{i}}).$ }
\end{itemize}

Details are explained in \cite{Cas}. In  \cite{Kar} the reader will
find examples of transversal structures with applications to
algebraic topology.

{\noindent}We think of $\C(x_{0},\cdots ,x_{n}) $ as the
differential graded subspace of
$\des\bigotimes_{i=1}^{n}\C(x_{i-1},x_{i})$ consisting of
transversal $n$-tuples $\C(x_{0},\cdots ,x_{n}).$ Thus any $0$ or
$1$-tuple is automatically transversal. For $n\geq 2$ any closed
$n$-tuple in $\displaystyle\bigotimes_{i=1}^{n}\C(x_{i-1},x_{i})$
is homologous to a closed transversal $n$-tuple, and any
subsequence of a transversal sequence is transversal.

{\noindent}Let $D^{1}=\{y\in\mathbb{R}\mid |y|\leq1\}.$ The operad
of little intervals $\I=\{\I(n)\}_{n\geq1}$ is defined as follows:
for each $n\geq1,$ $\I(n)$ is the topological space
\[\I(n)= \Big\{(T_{y_{1},r_{1}}, \ldots ,T_{y_{n},r_{n}}) \ \Big |
\begin{array}{c}
y_{i}, y_{j}\in D^{1},\ 0\leq r_{i} < 1 \ \mbox{such that if }\
i\neq j\
 \\
 \overline{\im(T_{y_{i},r_{i}})}\cap \overline{\im(T_{y_{j},r_{j}})}=
\emptyset \ \mbox{ for all} \ 1\leq
i,j \leq n  \\
\end{array} \Big \} ,\]
where $T_{a,r}\colon D^{1}\to D^{1}$ is the affine transformation
 $T_{a,r}(y)=ry+a,$ for $x\in D^{1}$ and $r\geq0$ small
 enough.

\begin{figure}[h]
\begin{center}
\includegraphics[height=1.6cm]{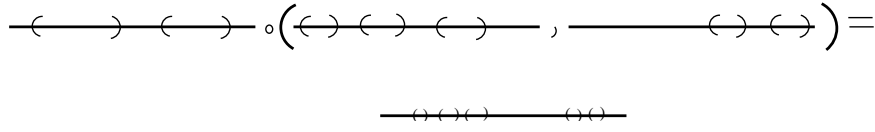}
\caption{\ Example of a composition $\gamma_{2}\colon \I(2)\otimes
\I(3)\otimes \I(2)\to \I(5) $. \label{figu1}}
\end{center}
\end{figure}

{\noindent} We regard $\I$ as a non-symmetric operad for the
purposes of this
 note. Notice that a configuration of little intervals is
 naturally numbered using the order in $D^{1}.$ Compositions in
 $\I=\{\I(n)\}_{n\geq1}$
are illustrated by Figure \ref{figu1}

{\noindent}We shall use the dg-vect operad of singular chains
$\Ca_{s}(\I)=\{\Ca_{s}(\I(n))\}_{n\geq 1}$ with compositions induced
by those on $\I.$
 A dg-precategory $\C$ is said to be a {\it
 transversal $1$-category} if a domain $\C^{\star}$ in
$\C$ is given together with maps $\theta_{n}\colon
\Ca_{s}(\I(n))\otimes \C(x_{0},\cdots ,x_{n}) \to \C(x_{0},x_{n})$
for any $x_{0}, \cdots ,x_{n} \in \ob( \C).$ This data should
satisfy the following axioms:
 1) $\theta_{1}(1\otimes(\ ))=\id\colon \C(x_{0},x_{1}))\to\C(x_{0},x_{1}),$
where $1$ denotes the identity in $\Ca(\I(1)).$ 2) For each
partition $n=n_{1}+ \cdots + n_{k},$ the map
$\displaystyle\bigotimes_{i=1}^{k}\Ca_{s}(\I(n_{i}))\otimes
(\displaystyle\bigotimes_{i=1}^{n} \C(x_{n_{i-1}},\cdots
,x_{n_{i}}))\to \displaystyle\bigotimes_{i=1}^{k} \C(x_{n_{i-1}}
,x_{n_{i}})$ obtained by including $\C(x_{0},\cdots ,x_{n})$ in
$\displaystyle\bigotimes_{i=1}^{k} \C(x_{n_{i-1}},\cdots
,x_{n_{i}}),$ shuffling, and applying $\theta ^{\otimes k}$ factors
through $\C(x_{n_{0}},x_{n_{1}},\cdots,x_{n_{k}}),$ see \cite{Cas}.
Using that $\Ho(\I(n))=k$ for all $n\geq1,$ it is shown in
\cite{Cas} that

\begin{thm}\label{Teo4}
{\em If $\C$ is a transversal $1$-category then $\Ho(\C)$ is a
category.}
\end{thm}

{\noindent}The concept of transversal $1$-algebra may be deduced
from that of a transversal category as follows: a dg-vect $A$ is a
transversal $1$-algebra if and only if the precategory $\C_{A}$
with only one object  $p,$ such that $\C_{A}(p,p)=A$ is a
$1$-category.

\begin{cor}\label{cor5}
 {\em If $A$ is a transversal $1$-algebra then
$\Ho(A)$ is an associative algebra.}
\end{cor}

{\noindent} Smooth manifolds with corners are defined in \cite{Mel}.
Given a smooth manifold $M$ we define the graded vector space
\[\Ca(M)=\bigoplus_{i\in \N} \Ca_{i}(M),\]
with $i\in \N,$ where $\Ca_{\it i}(M)$ denotes the complex
 vector space
\[ \frac{\Big\langle (K_{x},x)\colon \begin{array}{c}
K_{x}\ \mbox{is a compact connected oriented}\ i \mbox{-manifold with }   \\
\mbox{corners and}  \ x\colon K_{x} \to M \
\mbox{is a smoth map} \\
\end{array}
\Big\rangle}{\langle (K_{x},x)-(-K_{x},x) \rangle}.\]

{\noindent}For $K$ an oriented manifold, $-K$ denotes the same
manifold provided with the opposite orientation. We let $\partial
\colon \Ca_{i}(M) \to \Ca_{{i}-\rm {1}}(M)$ be the map given by
$\partial(K_{x},x)=\sum_{c}(\overline{c},x_{\overline{c}}),$ where
the sum ranges over the connected components  of $\partial_{1}
K_{x}$ provided with the induced boundary orientation and
$x_{\overline{c}}$ denotes the restriction of $x$ to the closure of
$c.$ In \cite{Cas} is proven that

\begin{thm}
 {\em $(\Ca(M),\partial)$ is a differential $\mathbb{Z}$-graded
$k$-vector space. Moreover}
\[\Ho(\Ca(M),\partial)=\Ho(M):=\mbox{\em singular homology of } \
M.\]
\end{thm}

\section{Homological quantum field theory}\label{sec2}

 This work was motivated by our desired to understand and
generalized string topology introduce by M.Chas and D. Sullivan
\cite{SCh}, see also \cite{CJ} and \cite{S}. Indeed the category
$\Co^{M}_{1,r}$ of homological one dimensional restricted cobordisms
enriched over a manifold $M$ constructed in Section \ref{sec3}
includes open string topology on $M$ as a particular case when one
considers objects of the form $c \colon [1] \longrightarrow D(M).$

The theory of cobordisms was introduced by Rene Thom \cite{RT}.
Using cobordisms Michael F. Atiyah \cite{MA} and G. B. Segal
\cite{Se} introduced the axioms for topological quantum field
theories TQFT and conformal field theories CFT, respectively. V.
Turaev \cite{TU1} introduced the axioms for homotopical quantum
field theories. Essentially the axioms of Atiyah may be summarized
as follows: 1) Consider the monoidal category $\Co_{n}$ of
$n$-dimensional cobordism. 2) Define the category of TQFT as the
category $\Mfu(\Co_{n},\ve)$ of monoidal functors $\f\colon
\Co_{n}\to \ve.$ The category of  Segal's  CFT may also be recast as
a category of monoidal functors.

{\noindent}We shall define the category $\HL$ of homological
quantum field theories as follows

{\noindent}1) Construct a $1$-category $\sC \OB_{n}^{M}$ for each
$n\geq1$ and compact oriented smooth manifold $M.$

{\noindent}2) Consider the associated category
$\Co_{n}^{M}=\Ho(\sC\OB_{n}^{M}).$

{\noindent}3) Define the category $\HL(M,n)$ to be the category of
monoidal functors $\Mfu(\Co_{n},\ve).$

\begin{figure}[ht]
\begin{center}
\includegraphics[height=2.1cm]{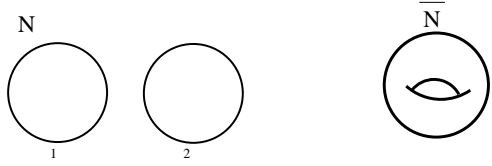}
\caption{ \ Manifold $N$ and its completion
$\overline{N}.$\label{fieti}}
\end{center}
\end{figure}

{\noindent}Let us first define the transversal $1$-category $\sC
\OB_{n}^{M}.$   By we denote $D(M)$ the set of all embedded
connected oriented submanifolds of $M.$ $\pi_{0}(M)$ denotes the
set of connected components of $M.$  We define the completion of a
manifold $M$ to be $\overline{M}=\des\prod_{c\in \pi_{0}(N)}c.$
Figure \ref{fieti} gives an example of a manifold and its
completion.

 {\noindent}Objects in $\sC \OB_{n}^{M}$ are triples $(N,f,<)$ where
 $N$ is a compact oriented  manifold of dimension $n-1,$
 $f\colon \pi_{0}(N) \to D(M)$  is any map, and $<$ is a linear ordering on
 $\pi_{0}(N).$ By convention the empty
set is assumed to be a $n$-dimensional manifold for all $n\in
\mathbb{N}.$

{\noindent}Let $(N_{i},f_{i},<_{i})$ and $(N_{o}, \ f_{o},<_{o})$
be objects in $\sC \OB_{n}^{M}.$ We set
\[ \sC \OB_{n}^{M}((N_{i},f_{i},<_{i}), \ (N_{o},f_{o},<_{o}))=\overline{\sC
\OB_{n}^{M}}\diagup \backsim\] where $\overline{\sC \OB_{n}^{M}}$
is the set of triples $(P,\alpha,\xi)$ such that

{\noindent}1) $P$ is a compact oriented smooth manifold with
corners of dimension $n.$

{\noindent}2) $\alpha \colon N_{i}\bigsqcup N_{o}\times [0,1]\to
\im (\alpha)\subseteq P$ is a diffeomorphism and $\alpha \mid_{
N_{i}\bigsqcup N_{o}}\to
\partial  P$ is a diffeomorphism such that $\alpha|_{N_{i}}$
reverses orientation and $\alpha|_{N_{o}}$ preserves orientation.

{\noindent}3) $\xi \in \Ca(\ma (P,M)_{f_{i}, f_{o}}),$ where
$\ma(P,M)_{f_{i}, f_{o}}$ is the set of all smooth maps $g\colon
P\to M$ such that for $c\in \pi_{0}(N_{j}), \ \xi(a, )$ is a
constant map with value in $f_{j}(c)$ on an open neighborhood of
$c,$ for $j=i,o.$  We considered $\Ca(\ma (P,M)_{f_{i}, f_{o}})$
with degrees shifted down by
  $\dim N_{0}.$ Thus $\xi$ is a smooth map $\xi\colon
K_{\xi}\times P\to M$ for some compact oriented manifold with
corners $K_{\xi}.$

{\noindent}If $\alpha\in \ma(PM)_{f_{i},f_{o}}$ we define
$e_{i}(\alpha)\in\des \prod_{c\in\pi_{0}(N_{i})}f_{i}(c)$ by
$e_{i}(\alpha)(c)=e_{i}(x)$ for $x\in c;$ similarly
$e_{o}(\alpha)\in\des \prod_{c\in\pi_{0}(N_{o})}f_{o}(c)$ is given
by $e_{o}(\alpha)(c)=e_{o}(x)$ for $x\in c.$ Triples
$(P_{1},\alpha_{1}, \xi_{1} )\sim (P_{2},\alpha_{2}, \xi_{2} )$ in
$\overline{\sC \OB_{n}^{M}}((N_{i},f_{i},<_{i}), \
(N_{o},f_{o},<_{o}))$ if and only if there is an orientation
preserving diffeomorphism $\varphi\colon P_{1} \to P_{2}$ such
that $\varphi \circ \alpha_{1}=\alpha_{2},$ and
$\varphi_{\star}(\xi_{1})=\xi_{2}.$

{\noindent}$\sC \OB_{n}^{M} ((N_{0},f_{0},<_{0}),\cdots ,
(N_{k},f_{k},<_{k}))$ is the set of all $k$-tuples
$\{(P_{i},\alpha_{i}, \xi_{i})\}_{i=1}^{k}$ such that
$(P_{i},\alpha_{i},\xi_{i})\in \sC \OB_{n}^{M}
((N_{i-1},f_{i-1},<_{i-1}), (N_{i},f_{i},<_{i}))$ and for all
$1\leq i \leq j \leq k$ the map
\[
\begin{array}{c}
  \txt \footnotesize {$\!\!\!\!\!\!\!\!\!\!\!\!\!\!\!\!\!\!
\!\!\!\!\!\!\!\!\!\!\!\!\!\!\!\!\!\!\!\!\!\!\!\!\!e (\xi_{1},
 \ldots , \xi_{n})\colon
 \des\prod_{i=1}^{k}K_{\xi_{i}}
\rightarrow
\des \prod_{i=1}^{k-1}\overline{N}_{i}\times\overline{ N}_{i}$}\\
\txt \footnotesize {$e (\xi_{1}, \ldots , \xi_{k}) (c_{1},\ldots ,
c_{k}) = (e_{o}(\xi_{1}(c_{1}))
,e_{i}(\xi_{2}(c_{2})),e_{o}(\xi_{2}(c_{2})),
\ldots ,e_{i}(\xi_{k}(c_{k}))) $}\\
\end{array}\]
is transversal to $\txt \footnotesize
{$\Omega_{k}=\des\prod_{i=1}^{k}
\Delta^{N_{i}}_{2}\subset\des\prod_{i=1}^{k-1}\overline{N}_{i}\times
\overline{N}_{i}$}$ where $\Delta^{N_{i}}_{2}=\{(a,a)\in
\overline{N}_{i}\times \overline{N}_{i} \}$ for $1\leq i \leq
k-1.$ Clearly
\begin{equation}\label{eva2}
\!\!\!\txt \footnotesize {$e^{-1}(\Omega_{k})= $}\Big\{\txt
\footnotesize {$(c_{1}, \ldots ,c_{k}) \in
\des\prod_{i=1}^{k}K_{\xi_{i}}$} \Big |
\begin{array}{c}
\txt \footnotesize
{$e_{o}(\xi_{i}(c_{i}))=e_{i}(\xi_{i+1}(c_{i+1}))$}
 \\
\txt \footnotesize {$1\leq i \leq k-1$}\\
\end{array}\Big \}
\end{equation}

{\noindent}Since $e$ is a smooth map and $e(\xi_{1},
 \ldots , \xi_{k})\pitchfork\Omega_{k}$ then
\[K_{\xi_{1}}\times_{\overline{N}_{1}}K_{\xi_{1}}\times_{\overline{N}_{2}} \cdots
\times_{\overline{N}_{k-1}} K_{x_{k}}=e^{-1}(\Omega_{k})\] is a
manifold with corners.

\begin{figure}[h]
\begin{minipage}[t]{0.48\linewidth}
\begin{center}
\includegraphics[height=2.7cm]{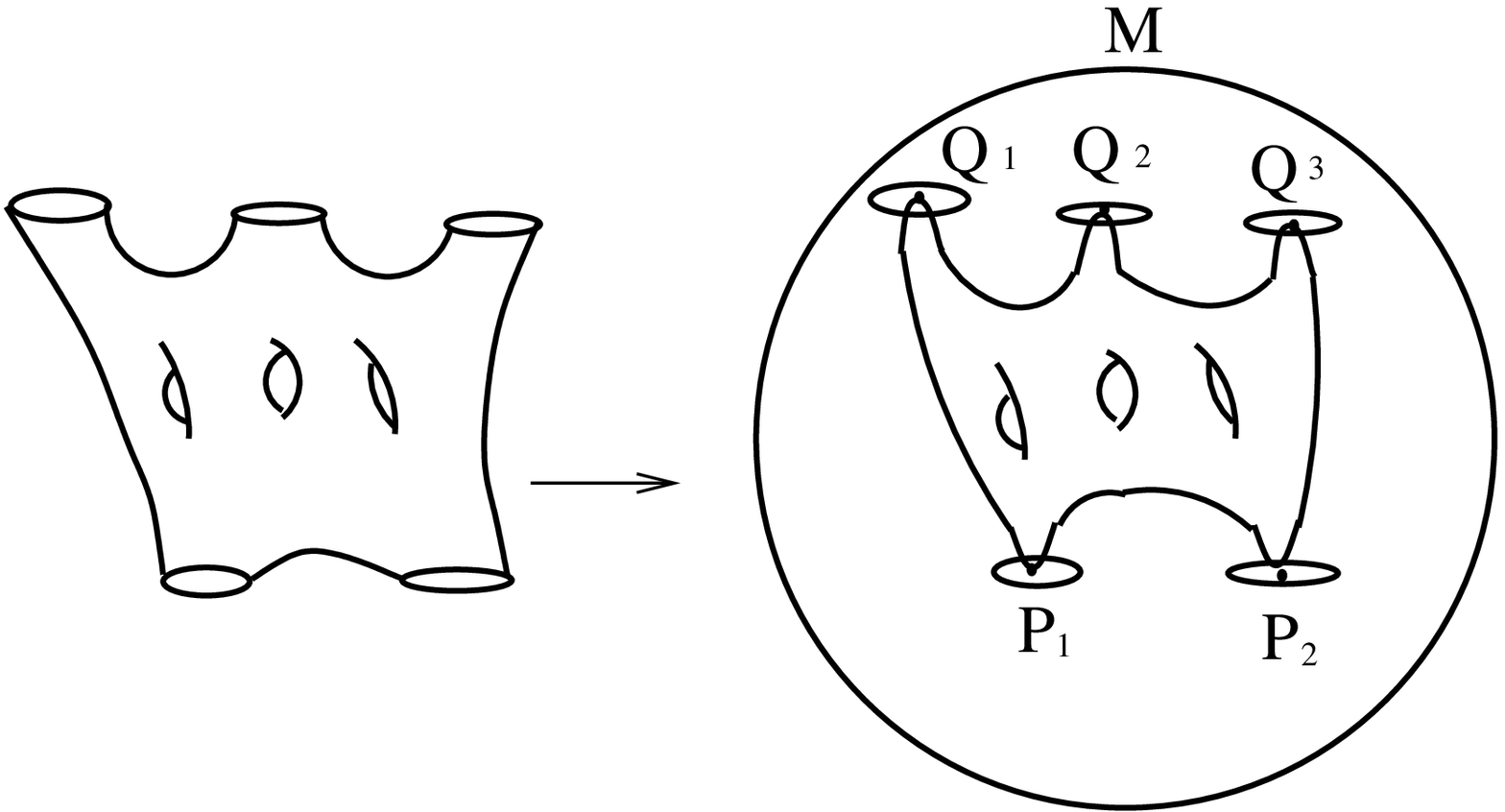}\!\!\!\!\!\!\caption{A
$n$-cobordism enriched over $M.$}\label{Bor3}
\end{center}
\end{minipage}
\begin{minipage}[t]{0.5\linewidth}
\begin{center}
\includegraphics[height=3.5cm]{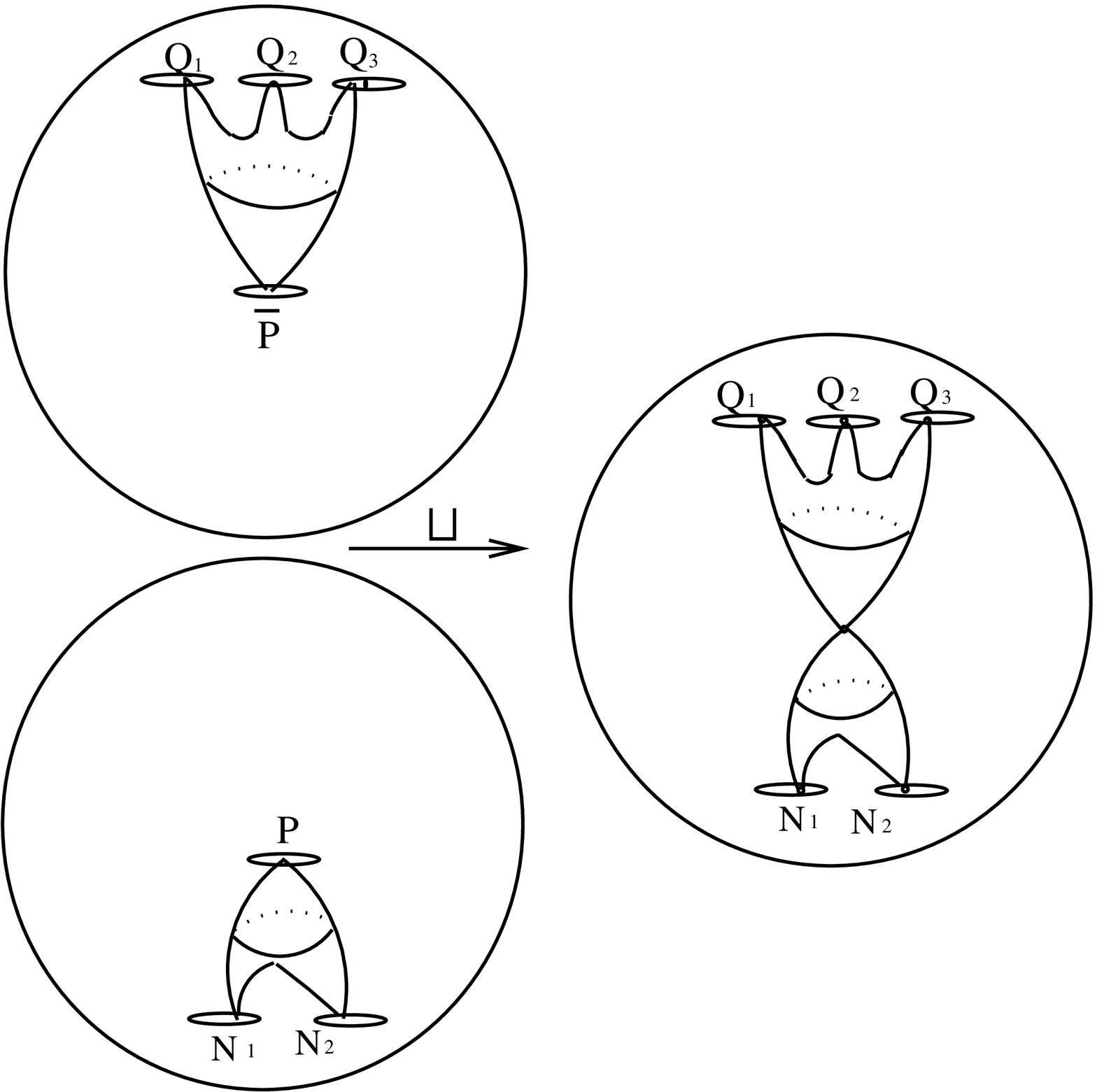}\!\!\!\!\!\! \!\!\!\!\!\!\ \caption{Composition of
$n$-cobordisms enriched over $M.$}\label{Bor4}
\end{center}
\end{minipage}
\end{figure}

{\noindent}Given $a\in \Ca_{s}(\I(k)),$
$(P_{i},\alpha_{i},\xi_{i})\in
 \sC \OB_{n}^{M}
((N_{i-1},f_{i-1},<_{i-1}), (N_{i}, f_{i}, <_{i}))$ for $1\leq
i\leq k,$ the composition morphism
$a((P_{1},\alpha_{1},\xi_{1}),\cdots
,(P_{k},\alpha_{k},\xi_{k}))\in \sC \OB_{n}^{M} ((N_{0},
f_{0},<_{0}), (N_{k},f_{k},<_{k})))$ is the triple
$(a(P_{1},\cdots , P_{k}),a(\alpha_{1},\cdots ,
\alpha_{k}),a(\xi_{1},\cdots , \xi_{k}))$ where
\[\begin{array}{c}
  a(P_{1},\cdots , P_{k})=P_{1}\des\bigsqcup_{N_{1}}
\cdots \des\bigsqcup_{N_{k-1}}P_{k} \\
  a(\alpha_{1},\cdots ,
\alpha_{k})=\alpha_{1}\mid_{N_{0}}\bigsqcup
\alpha_{k}\mid_{N_{k}} \\
          \\
  K_{a}(\xi_{1},\cdots , \xi_{k})=K_{a}\times
  K_{\xi_{1}}\des\times_{\overline{N}_{1}}K_{\xi_{2}}\times
\cdots \des\times_{\overline{N}_{k-1}}K_{\xi_{k}} \\
\end{array}\]
The map $a(\xi_{1},\cdots , \xi_{k})\colon K_{a}(\xi_{1},\cdots ,
\xi_{k})\times P_{1}\des\bigsqcup_{N_{1}}\cdots
\des\bigsqcup_{N_{k-1}}P_{k}\to M$ is given by $a(\xi_{1},\cdots ,
\newline \xi_{k})((s, t_{1}\cdots ,t_{k}),u)=\xi_{i}(t_{i})(u)$ for any
tuple $(s, t_{1}\cdots ,t_{k})$ in $K_{a}\times
K_{\xi_{1}}\times_{\overline{N}_{1}} \cdots
\times_{\overline{N}_{k-1}}  K_{\xi_{k}}$ and $u\in P_{i}.$
 Figure $\ref{Bor3}$ represents a $n$-cobordism
enriched over $M$ and Figure $\ref{Bor4}$ shows a composition of
$n$-cobordism enriches over $M.$

{\noindent}$\Co^{M}_{n,r}$  is the full subcategory of
$\Co^{M}_{n},$ such that the empty set is not included as an object.

\begin{prop}
 {\em $(\Co_{n}^{M}, \ \sqcup, \ \emptyset)$ is a monoidal category
with disjoint union $\sqcup$ as product and empty set as unit.
$ \ \Co^{M}_{n,r}$
 is a monoidal category without unit. }
\end{prop}

{\noindent}Given monoidal categories $\C$ and $\d$ we denote by
$\Mfu(\C,\d)$ the category of monoidal functors from $\C$ to $\d.$

\begin{defi}
{\em $\HL (M,d)=\Mfu(\Co_{n}^{M},\ve)$ and $\HL_{r} (M,d)=
\Mfu(\Co_{n.r}^{M},\ve).$ $\HL_{d}(M,d)$ is the category of}
homological quantum field theories of dimension $d.$ $\HL_{r}
(M,d)$ {\em is the category of} restricted homological quantum
field theories of
 dimension $d.$
\end{defi}

\vspace{0.3cm}

\section{One-dimensional \HL}\label{sec3}

\begin{figure}[ht]
\begin{center}
\includegraphics[height=4cm]{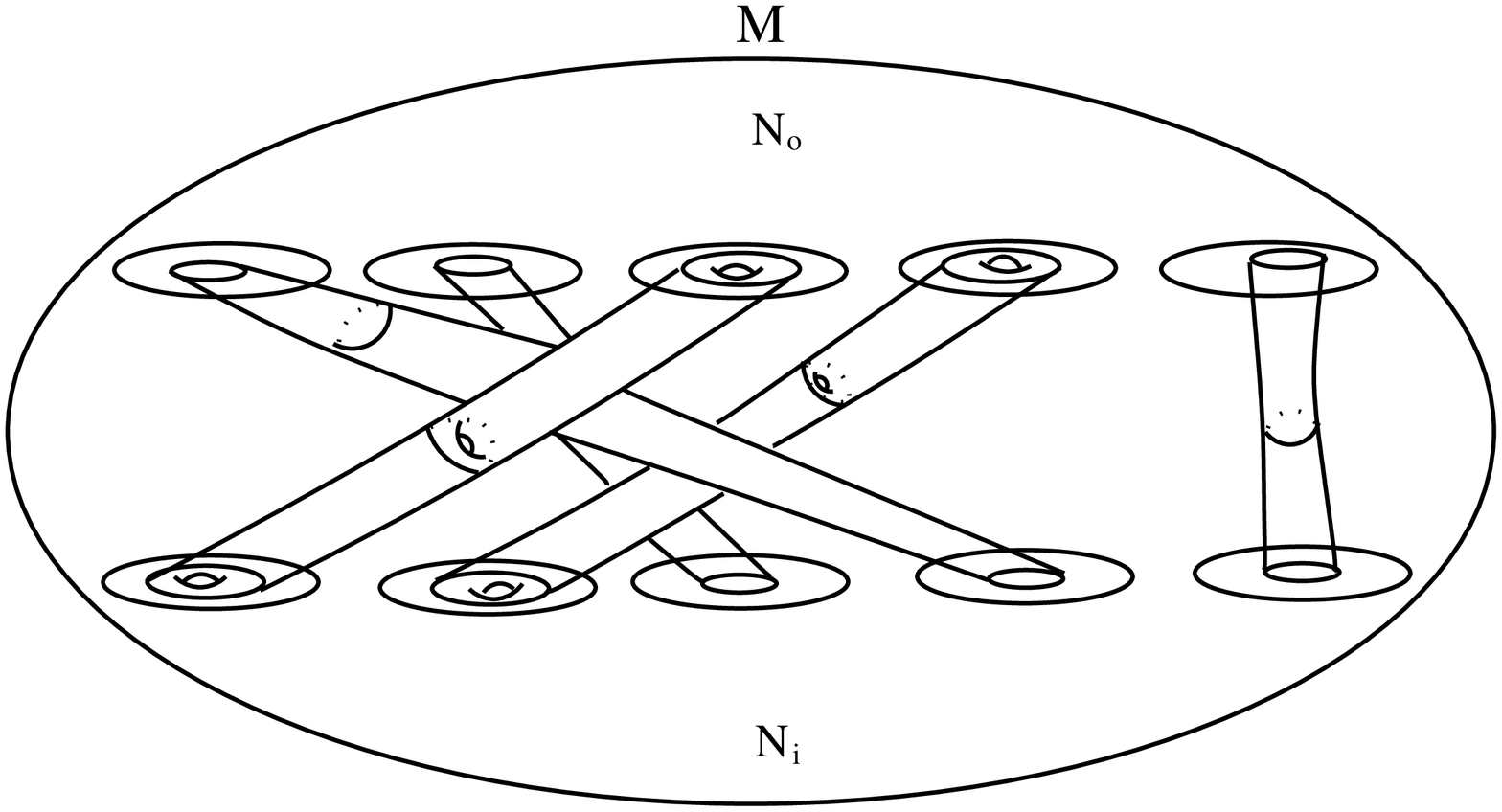}
\caption{ \ A morphism in $\Co_{1}$. \label{strin}}
\end{center}
\end{figure}

For $n\in \mathbb{N}$ we set $[n]=\{1,\cdots ,n\}.$ Let $M$ be a
compact oriented smooth manifold and $\pi\colon P\to M$ be a
principal fiber $G$-bundle for $G$ a compact Lie group. Let
$\scr{A}_{P}$ denote the space of planar connections on $P,$ and let
$\Lambda \in \scr{A}_{P}$  be a connection. If $\gamma \colon I \to
M$ is a smooth curve on $M$ and $x\in P$ is such that
$\pi(x)=\gamma(0),$ then we let $P_{\Lambda}(\gamma,x)$ be
$\widehat{\gamma}(1)$ where $\widehat{\gamma}$ is the horizontal
lift of $\gamma$ with respect to $A$ such that
$\widehat{\gamma}(0)=x.$ The goal of this section is to construct a
map $\Ho\colon \scr{A}_{P} \to \HL_{r}(M,1).$ For each connection
$\Lambda \in \scr{A}_{P}$ we construct a one-dimension restricted
homological quantum field theory $\Ho_{\Lambda}\colon
\Co_{1,r}^{M}\to \ve.$ For $c\in \ob(\Co^{M}_{1,r})$ we let
$\overline{c}$ be $\des \prod_{i\in[n]}c(i).$ Notice that an object
in $\Co^{M}_{1,r}$ is a map $c\colon [n]\to D(M)$ for some $n\in
\mathbb{N}^{+}.$ On objects $\Ho_{\Lambda}\colon
\ob(\Co_{1,r}^{M})\to \ob(\ve)$ is the map given by
\[\Ho_{\Lambda}(c)=\Ho(\des\prod_{i\in[n]} P\mid_{c(i)})\]
where $P\mid_{c(i)}$  denotes the restriction of $P$ to
$c(i)\subseteq M.$ Notice that a morphism in $\Co_{1}$ from $c$ to
$d$ consists of a pair $(\alpha, t)$ where $\alpha$ is a permutation
$\alpha \colon [n]\to [n],$ and $t=\{t_{1}\cdots t_{n}\}$ is such
that $t_{i}\in\Ho(\ma(\I_{c},M)_{c_{i},d(\alpha(i))})$ where $I_{c}$
is a $1$-dimensional compact manifold with boundary (a closed
interval). Figure \ref{strin} shows a picture of a morphism in
$\Co_{1}$.

{\noindent}The map
\[\Ho_{\Lambda}(\alpha, t)\colon \Ho(\des\prod_{i\in[n]} P\mid_{c(i)})\to
\Ho(\des\prod_{i\in[n]} P\mid_{d(i)}) \] is defined as follows:
consider the natural projection $\pi\colon\des\prod_{i\in[n]}
P\mid_{c(i)}\to \des\prod _{i\in[n]} c(i),$ and let $x$ be a chain
$x\colon K_{x}\to\des\prod_{i\in[n]} P\mid_{c(i)},$ where
$x=(x_{1},\dots ,x_{n}).$ The domain of $\Ho_{\Lambda}(\alpha,
t)(x)$ is given by
\[K_{\Ho_{\Lambda}(\alpha, t)(x)}=K_{x}\times_{\overline{c}}
\des\prod _{i\in[n]}K_{t_{i}}\] The map
$\Ho_{\Lambda}(\alpha,t)(x)\colon K_{\Ho_{\Lambda}(\alpha,t)(x)}\to
\des\prod_{i\in[n]} P\mid_{d(i)}$ is given by
\[[\Ho_{\Lambda}(\alpha, t)(x)](y;y_{1},\cdots
y_{n})=\des\prod_{i\in[n]} P(t_{i}(y_{i}),x_{i})\] where $y\in
K_{x}, \ y_{i}\in K_{t_{i}}.$

\begin{thm}\label{t6}
 {\em For any $\Lambda \in \scr{A}_{P}$ the map $\Ho_{\Lambda}$ defines on
 $1$-dimensional homological quantum field theory.}
\end{thm}
\begin{proof}
We need to show that $\Ho_{\Lambda}((\beta,s)\circ (\alpha,
t))=\Ho_{\Lambda}(\beta,s)\circ \Ho_{\Lambda}(\alpha,t).$ Since
\[K_{\Ho_{\Lambda}((\beta,s)\circ (\alpha, t))}(x)=
K_{x}\times_{\overline{N}_{i}}(\des\prod_{c\in\pi_{0}(N_{i})}K_{t(c)}\times_{\overline{N}_{m}}
\des\prod_{d\in\pi_{0}(N_{m})}K_{s(d)})
\]
and
\[K_{\Ho_{\Lambda}(\beta,s)\circ \Ho_{\Lambda}(\alpha,t)}=(K_{x}\times_{\overline{N}_{i}}
\des\prod_{c\in\pi_{0}(N_{i})}K_{t(c)})\times_{\overline{N}_{m}}\des\prod_{d\in\pi_{0}(N_{m})}K_{s(d)}\]
the domains $\Ho_{\Lambda}((\beta,s)\circ (\alpha, t))$ and
$\Ho_{\Lambda}(\beta,s)\circ \Ho_{\Lambda}(\alpha,t)$ agree. The
corresponding maps also agree since $P(\gamma \circ \beta,
x)=P(\alpha, P(\beta,x)).$
\end{proof}

\section{Homological matrices}\label{sec4}

\begin{figure}[ht]
\begin{center}
\includegraphics[height=2.3cm]{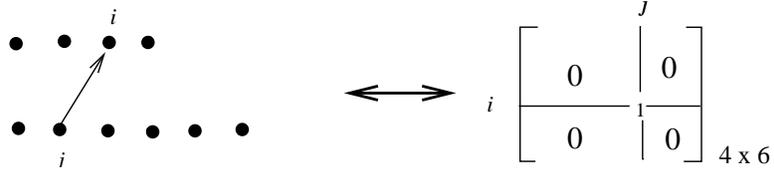} \caption{ \ Graphical representation
of matrices  \label{mat1}}
\end{center}
\end{figure}

In this section we develop a higher dimensional analogue of matrices
and matrix multiplication. D\'iaz and Pariguan in \cite{RE} used the
following graphical representation of matrices.  Let
$\textrm{Digraph}(m,n)$ be the $\mathbb{C}$-vector space generated
by the set of directed bipartite graphs $\Gamma$ with one edge
starting in $[m]$ and ending in $[n],$ so $\Gamma$ is given by an
element $\Gamma \in [m]\times [n].$ $\textrm{Digraph}(m,n)$ and the
vector space $M_{m\times n}(\mathbb{C})$ of complex matrices of
format $m\times n$ are naturally isomorphic via the correspondence
shown in Figure \ref{mat1}.

{\noindent}Fix a compact oriented manifold $M.$ We proceed to
construct a higher dimensional homological analogue of matrices by
generalizing the spaces $\textrm{Digraph}(m,n)$ as follows: 1) We
replace each point of  the set $[m]$ by oriented embedded
submanifold of $M.$ So the indices of a generalized matrix are maps
$c\colon [m]\to D(M).$ 2) The vector space $\textrm{Digraph}(m,n)$
is replaced by the space $\hg(c,d)$ of homological bipartite graphs
with one packet of edges,  where $c\colon [m]\to D(M)$ and $d\colon
[n]\to D(M).$ The space $\hg(c,d)$ is given by
\[\hg(c,d)=\des\bigoplus_{j\in[m] i\in[n]}\Ho(\ma(I,M)_{c(j),d(i)})\]
with degrees in $\Ho(\ma(I,M)_{c(j),d(i)})$ shifted by
$\dim(d(i)).$ 3) Given $A\in \hg(d,e)$ and $B\in \hg(c,d)$ then we
define the product $AB\in \hg(c,e)$ of $A$ and $B$ by the rule
\[(AB)_{ki}=\des\sum_{j\in[m]}A_{kj}\circ B_{ji}\]
for all $i\in[n], j\in[m] \mbox{ and } k \in [p]$ where $B_{ji}\in
\Ho(\ma(I,M)_{c(i),d(j))}, \ A_{kj}\in \Ho(\ma(I,M)_{d(j),e(k))}$
and the composition maps
\[\circ\colon \Ho (\ma(I,M)_{d(j),e(k)})\otimes \Ho (\ma(I,M)_{c(i),d(j)})\to
\Ho (\ma(I,M)_{c(i),e(k)})\] are defined as for $1$-dimensional
homological quantum field theories.

\begin{figure}[ht]
\begin{center}
\includegraphics[height=2.3cm]{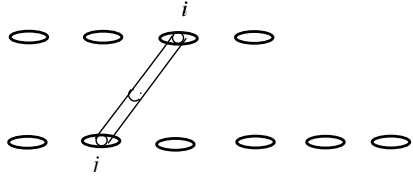} \caption{ \ Example of an element in $ \hg(c,d).$
  \label{mat2}}
\end{center}
\end{figure}

\begin{figure}[ht]
\begin{center}
\includegraphics[angle=90, height=2.7cm]{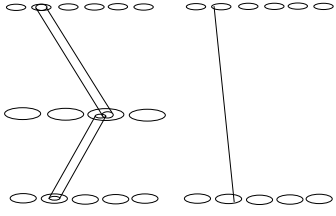} \caption{ \
Homological analog of matrix products \label{mat3}}
\end{center}
\end{figure}

{\noindent}For a map $c\colon[n]\to D(M),$ the space $\hg(c,c)$ is
an algebra which may be regarded as higher dimensional homological
analogue of $M_{n}(\mathbb{C}).$ A representation of $\hg(c,c)$ on
a $\mathbb{C}$-vector space $V$ is a linear map $\rho \colon
\hg(c,c)\to \en(V)$ such that $\rho(AB)=\rho(A)\rho(B).$ Figure
\ref{mat3} shows an example of product of homological matrix in
$\hg(d,e)$ with a homological matrix $\hg(c,d).$

\begin{thm}
 {\em $\des\bigoplus_{i \in [n]}\Ho(c(i)),$ where $\Ho(c(i))$ is shifted down by $\dim(c(i)),$
  is a representation of $\hg(c,c)$.}
\end{thm}
\begin{proof}
Let $A\in \des \bigoplus_{i,j\in [n]} \Ho (\ma(I,M)_{c(i),c(j)})$
and $e\in \des\bigotimes_{i \in[n]} \Ho (c(i))$ then we set
$(Ae)_{i}=\des\sum_{j\in[n]} A_{ij}\circ e_{j},$ where $A_{i,j}\in
\hg(c(j),c(i))$ and $e_{j}\in \Ho(c(j))$ and  the composition maps
 \[\circ \colon  \Ho
(\ma(I,M)_{c(i),c(j)})\otimes\Ho(c(j))\to \Ho(c(i))\] are define
as in Theorem \ref{t6}.
\end{proof}

{\noindent}In \cite{RE1} the $m$-symmetric power $\Si^{m}\al$ of
an algebra  $\al$ is defined as follows: as a vector space
$\Si^{m}\al=(\al^{\otimes m})_{S_{m}},$ ($S_{m}$-coinvariants.)
The product in $\Si^{m}\al$ is that as if $\overline{a_{1}\otimes
\cdots \otimes a_{m}}\in\Si^{m}\al$ and $\overline{b_{1}\otimes
\cdots \otimes b_{m}}\in\Si^{m}\al$ then
\[(\overline{a_{1}\otimes \cdots \otimes
a_{m}}) \cdot (\overline{b_{1}\otimes \cdots \otimes b_{m}})=
\frac{1}{m !} \des\sum_{\sigma \in S_{m}}
 {\rm
sgn}(a,b,\sigma)\overline{(a_{1}b_{\sigma^{-1}(1)})\otimes \cdots
\otimes (a_{m}b_{\sigma^{-1}(m)})}\] where ${\rm
sgn}(a,b,\sigma)=(-1)^{e}$ and $e=e(a,b,\sigma)= \sum_{i>j}|a_{i}
b_{\sigma^{-1}(j)}| + \sum_{\sigma(i)>\sigma(j)}|b_{i} b{_j}|.$

{\noindent}As shown in \cite{RE1} if $\al$ is the algebra
$\en(\R^{p\mid q})$ of supermatrices of dimension $(p\mid q)\times
(p\mid q)$ then $\Si^{m}(\en)(\R^{p\mid q})$ is the Schur
superalgebra (see \cite{Gre}) which controls the polynomial
representations of $\en(\R^{p\mid q}).$ We proceed to define  Schur
algebras in the homological context.

\begin{defi}
{\em For $c\colon [n]\to D(M)$ we set
$\Sc_{m}(c,c)=\Si^{m}(\hg(c,c)).$ }
 \end{defi}

\begin{figure}[ht]
\begin{center}
\includegraphics[height=1.5cm]{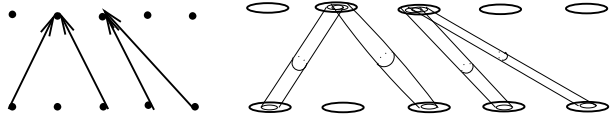} \caption{ \ Higher
dimensional analog of Schur algebras. \label{mat4}}
\end{center}
\end{figure}

\begin{figure}[ht]
\begin{center}
\includegraphics[height=3.6cm]{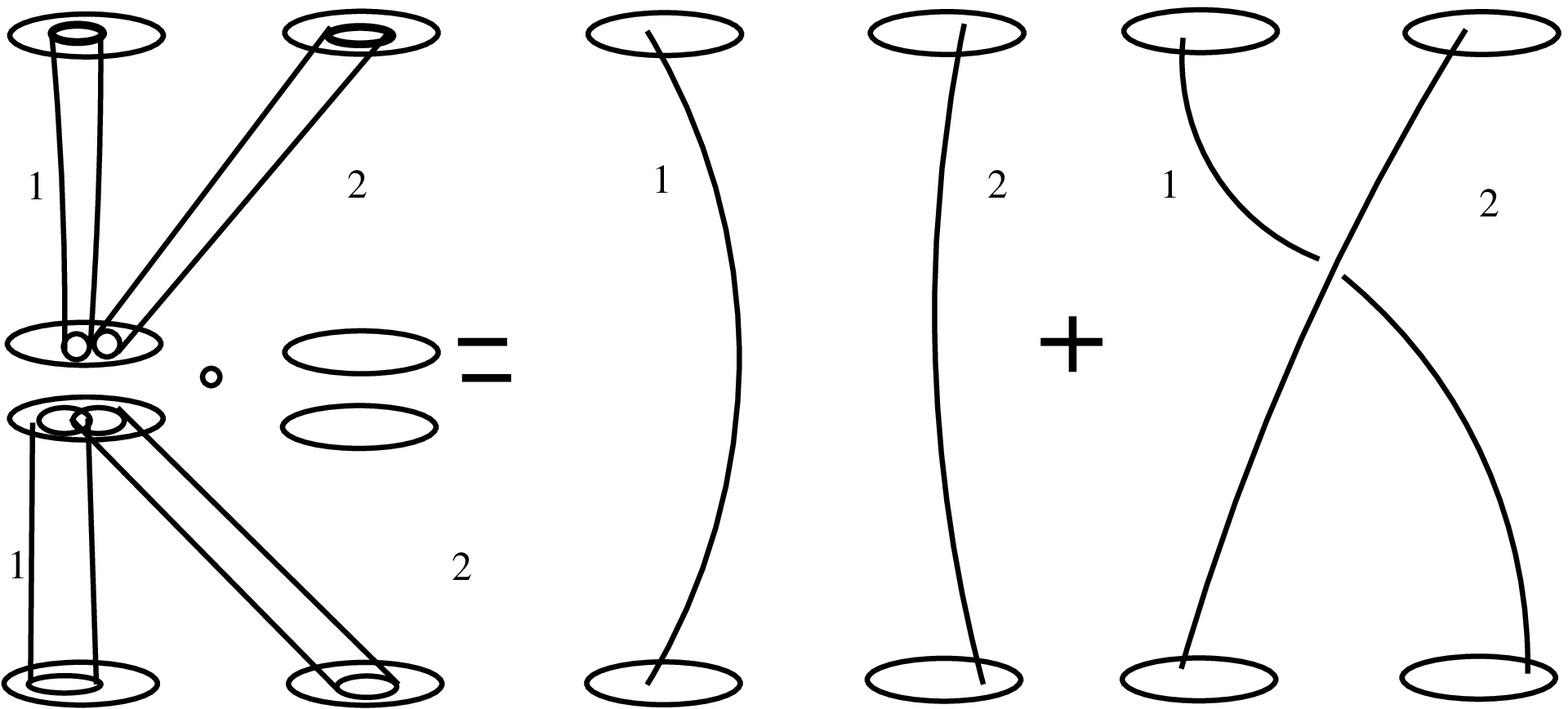} \caption{ \ Example of a product in
$\Si^{2}(\hg(1,1)).$ \label{matt3}}
\end{center}
\end{figure}

{\noindent}Figure \ref{mat4} shows an element  of $\Sc_{4}(c,c)=
\Si^{4}(\Ho(c,c)).$ Figure \ref{matt3} shows a product in
$\Sc_{2}(c,c).$

\begin{thm}
{\em If $V$ is a representation of $\Ho(c,c)$ then $\Si^{m}V$ is a
representation of
 $\Sc_{m}(c,c)$.}
\end{thm}
\begin{proof}
$\overline{a_{1}\otimes \cdots \otimes a_{m}}\in\Sc_{m}(c,c)$ and
$\overline{v_{1}\otimes \cdots \otimes v_{m}}\in\Si^{m} V$ then
\[\overline{a_{1}\otimes \cdots \otimes
a_{m}} \cdot \Big(\overline{v_{1}\otimes \cdots \otimes
v_{m}}\Big)= \frac{1}{m !} \des\sum_{\sigma \in S_{m}}
 {\rm
sgn}(a,v,\sigma)\overline{a_{1}(v_{\sigma^{-1}(1)})\otimes \cdots
\otimes a_{m}(v_{\sigma^{-1}(m)})}\] where ${\rm
sgn}(a,v,\sigma)=(-1)^{e}$ and $e=e(a,v,\sigma)= \sum_{i>j}|a_{i}(
v_{\sigma^{-1}(j)})| + \sum_{\sigma(i)>\sigma(j)}|v_{i} v{_j}|.$
\end{proof}

For the rest of the paper we work with \emph{even} dimensional
homology groups. Also we assume that the map $c\colon [n]\to D(M)$
is such that $\dim{(c(i))}$ is even for all $i\in[n].$

\begin{thm}
 {\em Let $c\colon [n]\to D(N)$ and  $C=\des\bigsqcup_{i\in[n]}c(i).$ There  is an algebra inclusion
  $i\colon \Co_{1,r}^{M}(C,C)\to \Sc_{m}(c,c).$ Thus any representation
 of $\Sc_{m}(c,c)$ induces a representation of $\Co_{1,r}^{M}(C,C)$.}
\end{thm}
\begin{proof}
Al element of $\Co_{1,r}^{M}(C,C)$ is a sequence
$a=(a_{i})_{i=1}^{n}$ where  $a_{i}\in \Ho(\ma(I,
M)_{c(i),c(\alpha(i))})$ for some $\alpha\colon [n]\to[n]$
biyective. The map $i$  is such that $i(a)=\overline{a_{1}\otimes
\cdots a_{m}}.$
\end{proof}

\begin{cor}
{\em $\Si^{m}(\des \bigoplus_{i=1}^{n}\Ho(c(i)))$ is a
representation of $\Co_{1,r}^{M}(C,C).$ }
\end{cor}

 \vspace{0.3cm}
\subsection*{Acknowledgment} We thank Delia Flores de Chela, Raymundo Popper, Eddy Pariguan, Sylvie
Paycha and Bernardo Uribe.

\bibliographystyle{amsalpha}

\end{document}